\documentclass[12pt]{ article}
\usepackage{amsmath}
\usepackage{amsthm}
\usepackage{amscd}
\usepackage{amssymb}
\usepackage{graphicx}
\newtheorem{thm}{Theorem}
\newtheorem{lem}{Lemma}
\newtheorem{cor}{Corollary}
\begin{document}
\title{ Khovanov Homology and the Twist Number of Alternating Knots}
\author{ Robert G. Todd\\ University of Iowa}
\maketitle

\subsection{Abstract}
In \cite{dasb} it was  shown that the sum of the absolute value of the second and penultimate coefficient of the Jones polynomial of an alternating knot is equal to the twist number of the knot. Here we give a new proof of  their result using a variant of Khovanov's homology that was defined by O. Viro for the Kauffman bracket. The proof is by induction on the number of crossings using the long exact sequence in Khovanov homology corresponding to the Kauffman bracket skein relation.

\section{ Categorification}

A categorification of the Kauffman Bracket is given in \cite{viro}. Given a diagram D with n crossings that have been ordered define a state to be an n-tuple whose entries are either 0 or 1. To each state associate a collection of circles by smoothing each crossing positively whose corresponding entry is a 0 and smoothing each crossing negatively whose corresponding entry is 1. Given a state s we can define an enhanced state S by assigning an orientation to each circle in the state. Let $|S|$ be the number of circles that occur after smoothing each crossing. Consider the following parameters.
\begin{equation} \sigma(S)= \sharp\{positive\ smoothings \} - \sharp\{negative\ smoothings\} \end{equation}
\begin{equation}
\tau(S)=\sharp \{clockwise\ oriented\ circles\} -\sharp\{counterclockwise\ oriented\ circles\}
\end{equation}

\begin{equation}
J(S)=\sigma(S) +2 \tau(S)
\end{equation}

 Define $C_{j,k}$ to be  a $\mathbb{Z}$-module whose generators are those enhanced states S with $j=\sigma(S)$ and $k=J(S)$.  Let $C_{k}=\oplus_{j}C_{j,k}$. Enhanced states $S_{1}$ and $S_{2}$ have non-zero incidence number only if the states $s_{1}$ and $s_{2}$ from which they result differ at a single crossing. Then we can obtain $S_{2}$ from $S_{1}$ by either joining two circles together or by splitting a circle apart. The incidence number of two enhanced states is zero unless the circles not involved in the change of soothing have the same orientation. Given the action and the orientation of the the circles of enhanced state $S_{1}$, the following rules describe the orientation of the new circles of enhanced state $S_{2}$;

1)  If  two counterclockwise oriented circles become one then the new circle should be oriented counterclockwise.

2) If two circles with opposite orientation become one then  the new circle 
should be clockwise

3)If one circle with clockwise orientation becomes two circles they  should have clockwise orientation.

4) If one circle with counterclockwise orientation becomes two circles they should have opposite orientation.

For each of the above cases, the incidence number $(S_{1};S_{2})=\pm1$.
To determine the sign, look to the order of the crossings. If the states for $S_{1}$ and $S_{2}$ differ in the $i^{th}$ entry let t be the number of 1's that appear before the $i^{th}$ entry. Then the sign is $(-1)^{t}$.

Khovanov originally defined these  chain groups  and maps with a different indexing and showed that the homology of this complex is an invariant of the knot diagram. The indexing given here was developed by Viro and gives a framed invariant. In fact we get the following formulation of the Kauffman bracket, up to a factor of $i=(-1)^{\frac{1}{2}}$.

\begin{equation}
<D>=\sum_{j ,k}(i)^{j} rank(H_{j,k}) A^{k}
\end{equation}

\section{A Short Exact Sequence}
In \cite{viro} there is a short exact sequence coming from the skein relation $$<D>=A <D_{-}> + A^{-1}<D_{+}>$$It is categorified by the following relation;
\begin{equation}
\alpha:C_{i,j}(D_{-})\longrightarrow C_{i-1,j-1}(D)
\end{equation} 
is the map which takes enhanced Kauffman states of $D_{-}$ to the enhanced Kauffman state of D whose circles and orientations are the same. Define another map
\begin{equation}
\beta: C_{i,j}(D) \longrightarrow C_{i-1,j-1}(D_{+})
\end{equation}
This map sends any state which has a negative crossing at c  to 0 while sending other states in C(D) to the states in $C(D_{+})$ which has corresponding circles and orientations of those circles.
Provided that the crossing which we have chosen to smooth is the the last one in the ordering the maps $\alpha$ and $\beta$ commute with the boundary operators making them into  homomorphisms of complexes. Thus $\alpha$ and $\beta$ form a short exact sequence of complexes
\begin{equation}
0\rightarrow C_{*,*}(D_{-}) \longrightarrow C_{*-1,*-1}(D) \longrightarrow C_{*-2,*-2}(D_{+})\rightarrow 0
\end{equation}

 As such we have a series of Long exact sequences of homology

\begin{equation} \stackrel{\partial}{ \longrightarrow} H_{i+1.j+1}(D_{-}) \longrightarrow H_{i,j}(D) \longrightarrow H_{i-1,j-1}(D_{+})\stackrel{\partial}{\longrightarrow}H_{i-1,j+1}(D_{-})\longrightarrow \ldots
\end{equation}

This short exact sequence can also be understood as coming from a shifted mapping cone of the map $$\beta : D_{+} \rightarrow D_{-}$$ which changes a positive crossing to a negative crossing. That is $C(D) =M(\beta)$. At each level,$\beta$ is either joining two circles or splitting a circle in two. In either case the polynomial degree of $\beta$ is 2.

\section{The Kauffman Bracket}
In \cite{lick}  a knot was defined to be plus adequate if the state with all 0 smoothings has more circles than any state with only one 1 smoothing, and a knot is minus adequate if the state with all 1 smoothings has more circles than any state with only one 0 smoothing. A knot or link is adequate if it is both plus and minus adequate. 
Assume that D is an alternating adequate link.
Having arrived at $D_{-}$ and $D_{+}$ from D as described, the Kauffman bracket of these three links can be written as follows
\begin{eqnarray}
<D_{-}>=a_{k-3}A^{k-3} +a_{k-7}A^{\mathit{j}_7} +\ldots + a_{l+1}A^{l+1} \\
<D>=a_{k}A^{k}+a_{k-4}A^{k-4}+\ldots +a_{l+4}A^{l+4}+a_{l}A^{l}\\
<D_{+}>=a_{k-1}A^{k-1} + a_{k-5}A^{k-5}+\ldots+a_{l+3} A^{l+3}
\end{eqnarray}

Equation (4) shows that we may express the coefficients of the Kauffman bracket in the following manner
\begin{equation}
a_{j}=\sum_k (i)^{k} rank(H_{k,j})
\end{equation}

In the Homology group $H_{k,j}$ call the first index  the homological degree and the second index the polynomial degree.

\section{G and its Dual $G^{*}$}
Let L be any alternating link. Then to that link we may associate a checkerboard shading by coloring regions such  that we sweep out the black region when rotating the over crossing counter clockwise. 
\begin{center}
\begin{picture}(29,29)
\includegraphics{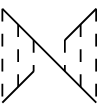}
\end{picture}
\end{center}

Then to such a coloring we can associate two one dimensional CW-complexes $G=(E,V)$ and $G^{*}=(E^{*},V^{*})$. Each vertex of G corresponds to a black region as a vertex of $G^{*}$ corresponds to a white region. An edge in G corresponds to a crossing which forms a corner between two black regions and an edge in $G^{*}$ corresponds to a crossing that forms a corner between two white regions. For example
\begin{center}
\begin{picture}(203,51)
\includegraphics{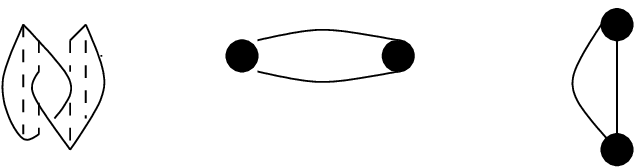}
\end{picture}
\end{center}

Define to edges to be parallel it they connect the same vertices. Then we define $\tilde{G}$, the reduced graph of G, to be the graph whose vertices are the same as those of G and whose edges correspond to equivalence classes of parallel edges in G. 

Now if $|E|$ is the number of edges and $|V|$ is the number of vertices. Let 
\begin{equation}
\psi(G)=|E| -|V| +1
\end{equation}

\section{The Geometry}
\begin{thm}
Let L be an alternating link. If L is plus adequate then
\begin{equation}
|a_{\mathit{k}-4}-a_{\mathit{k}}|=\psi(\tilde{G^{*}}) 
\end{equation}
And if L is minus-adequate Then
\begin{equation}
|a_{\mathit{l}+4} -a_{\mathit{l}}|=\psi(\tilde{G})
\end{equation}
\end{thm}

The proof will be by induction on the number of crossings of a knot diagram and will be based on the following lemma.

\begin{lem}

\begin{eqnarray}
i a_{k}(D) & = &a_{k-1}(D_{+}) \\
i a_{k-4}(D) & = & a_{k-5}(D_{+}) -a_{k-3}(D_{-})\\
i  a_{l}(D)&=& -a_{l+1}(D_{-})\\
i  a_{l+4}(D)&=& a_{l+3}(D_{+}) -a_{l+5}(D_{-})
 \end{eqnarray}  
 \end{lem}
 \begin{proof}[Proof of Lemma]
 Consider the long exact sequence
 \begin{equation}
 \stackrel{\partial}{\longrightarrow}H_{j+1,k+1}(D_{-}) \longrightarrow H_{j,k}(D) \longrightarrow H_{j-1,k-1}(D_{+})\stackrel{\partial}{\longrightarrow} H_{j-1,k+1}(D_{-}) \longrightarrow \ldots
 \end{equation}
 
 	Since D is adequate,in \cite{lick} it was shown that if we look at the number of circles created in the state of D where one crossing has been smoothed negatively and the others positively, it has one less circle than in the state of D where all crossing have been smoothed positively.  The state of $D_{-}$ where all crossings are smoothed positively corresponds to a state of D where one crossing has been smoothed negatively. Thus there is one less circle in the state of $D_{-}$ with all crossings smoothed positively than in the state of D with all positive smoothings. Recall that a particular enhanced state S  of a diagram is a generator of the Homology group with polynomial degree $\sigma(S) + 2\tau(S)$. Thus the highest possible polynomial degree occurs when all the crossings are smoothed positively and all the circles are oriented clockwise. Since the highest polynomial degree of D is k  and the state which contributes to the homology group of $D_{-}$ with highest polynomial degree has one less circle, this highest polynomial degree must be less than k.
	The above argument shows that $H_{*,k+1}(D_{-})$ is zero. Since D is adequate the same argument shows that the only homology group that contributes to the coefficient $a_{j}(D)$ is $H_{i,j}(D)$.  We are left with 

\begin{equation}
0\longrightarrow H_{j,k}(D) \longrightarrow H_{j-1,k-1}(D_{+}) \longrightarrow 0
\end{equation}

 Without loss of generality suppose that $(i)^{j-1}=-1$ . Thus

 \begin{eqnarray}
 rank(H_{j,k}(D))&=&(-i)a_{k}(D)\\
 rank(H_{j-1,k-1}(D_{+}))&=&a_{k-1}(D_{+})\\
  \end{eqnarray}
 The above  exact sequences give (16). 

Consider the long exact sequence of homology

\begin{equation}  
\begin{CD}
  @>\partial>>H_{j+1,k-3}(D_{-})@>\alpha_{*}>>H_{j,k-4}(D)@>\beta_{*}>>H_{j-1,k-5}(D_{+})\\@>\partial>>H_{j-1,k-3}(D_{-})@>\alpha_{*}>>H_{j-2,k-4}(D)@>\beta_{*}>>H_{j-3,k-5}(D_{+})@>\partial>>\ldots 
 \end{CD}
\end{equation}

Again, the alternating sum of the ranks of the homology groups being zero gives

\begin{multline}
(\ldots-rank(H_{j+1,k-3}(D_{-}))+rank(H_{j-1,k-3}(D_{-})-\ldots)\\
+\\(\dots +rank(H_{j,k-4}(D))-rank(H_{j-2,k-4}(D))+\ldots)\\
+\\(\ldots-rank(H_{j-1,k-5}(D_{+}))+rank(H_{j-3,k-5}(D_{+}))-\ldots)=0
\end{multline}

With the same sign conventions as above we recognize the following coefficients of the Kauffman bracket 

\begin{equation}
a_{k-3}(D_{-}) + i a_{k-4}(D) - a_{k-5}(D_{+})=0
\end{equation}

This gives us equation (17). Similar arguments give equations (18) and (19).

\end{proof}

\begin{proof}[Proof of Theorem]
Notice that the Hopf link is adequate, $\psi(\tilde{G})=0$, and $\psi(\tilde{G^{*}})=0$. $<Hopf>=-A^{6}-A^{2}-A^{-2}-A^{_6}$. Consequently,
\begin{eqnarray}
|a_{2}-a_{6}|&=&0=\psi(\tilde{G}^{*})\\
|a_{-2}-a_{-6}|&=&0=\psi(\tilde{G})
\end{eqnarray}
Suppose  the theorem holds for  alternating links with less than n crossings, and that L is an reduced alternating n crossing link.

 Two crossings are {\em equivalent} if there is a bi-gon connecting them. Below is an example of equivalent crossings. 
\begin{center}
 \begin{picture}(40,40)
\includegraphics{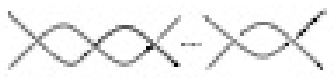}
\end{picture}
 \end{center}
 
 The number of such equivalence classes defines the twist number of a link. 

 The  crossing  to be smoothed in L is v, as indicated below.
\begin{center}
\begin{picture}(47,47)
\includegraphics{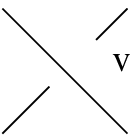}
\end{picture}
\end{center}

Assume  that v constitutes an equivalence class by itself. That is, there is no bi-gon connecting it to any other crossing. Then both $D_{-}$ and $D_{+}$ are alternating. Since there are no nugatory crossings in D then we can only create a nugatory crossing after smoothing at v if we created a kink. In order to create a kink v would have to have been equivalent to some other crossing which it is not. Thus $D_{-}$ and $D_{+}$ are reduced. 

In \cite{lick} it was shown that any alternating link with no nugatory crossings is adequate. Call a diagram with no nugatory crossings reduced. Both $D_{-}$ and $D_{+}$ are alternating and reduced and so they are adequate. Thus the only homology group which contributes to $a_{k-1}(D_{+})$ is $ H_{j-1,k-1}(D_{+})$, that is $ a_{k-1}(D_{+})=(i)^{j-1}$. Similarly for $ D_{-}$ the only Homology group that contributes to $a_{l+1}(D_{-})$ is $ H_{m+1,l+1}(D_{-})$
so $ a_{l+1}(D_{-})=(i)^{m+1}$. From lemma I we have 

\begin{eqnarray}
|a_{k-4}(D)-a_{k}(D)| & = & |a_{k-5}(D_{+})-a_{k-3}(D_{-})-a_{k-1}(D_{+})| \\
|a_{l+4}(D)-a_{l}(D)| & = & |a_{l+3}(D_{+})-a_{l+5}(D_{-})+a_{l+1}(D_{-})| 
\end{eqnarray}

Since $D_{-}$ is adequate the only homology group that contributes to $a_{k-3}(D_{-})$ is $H_{j-1,k-3}(D_{-})$. So $a_{k-3}(D_{-})=(i)^{j-1}=a_{k-1}(D_{+})$. When comparing the coefficients of the Jones polynomial and the Kauffman bracket in terms of highest and lowest degree we can ignore the contribution to the Jones polynomial by the writhe and look only at the normalization of the Kauffman bracket by the factor $-A^{2}-A^{-2}$. The first coefficient is $a_{k-1}(D_{+})$ and the second is $a_{k-5}(D_{+})-a_{k-1}(D_{+})$. Because the Jones polynomial of an alternating link has alternating signs, the sign of $a_{k-5}(D_{+})-a_{k-1}(D_{+})$ is opposite the sign of $a_{k-3}(D_{-})$. The same argument can be made for the sign of $a_{l+5}(D_{-})-a_{l+1}(D_{-})$ and the sign of $a_{l+3}(D_{+})$. Thus,
\begin{eqnarray}
|a_{k-4}(D)-a_{k}(D)| & = & |a_{k-5}(D_{+})-a_{k-1}(D_{+})|+1 \\
|a_{l+4}(D)-a_{l}(D)|& = & | a_{l+5}(D_{-})-a_{l+1}(D_{-})| +1
\end{eqnarray}

Since $D_{+}$ and $D_{-}$ are reduced alternating links with n-1 crossings we can say that

\begin{eqnarray}
|a_{k-4}(D)-a_{k}(D)| & = &\psi(\tilde{G^{*}_{+}}) +1\\
|a_{l+4}(D)-a_{l}(D)|& = &\psi(\tilde{G_{-}}) +1
\end{eqnarray}

Smoothing a crossing like v positively ,which is its own equivalence class, deletes an edge in $\tilde{G^{*}}$ and smoothing it negatively deletes an edge in 
$\tilde{G}$. Thus
\begin{eqnarray}
|a_{j-4}(D)-a_{j}(D)| & = &\psi(\tilde{G^{*}_{+}}) +1=\psi(\tilde{G^{*}})\\
|a_{l+4}(D)-a_{l}(D)|& = &\psi(\tilde{G}_{-}) +1=\psi(\tilde{G})
\end{eqnarray}

 Suppose v is one crossing in an equivalence class of k crossings. 
  \begin{center}
 \begin{picture}(125,127)
\includegraphics{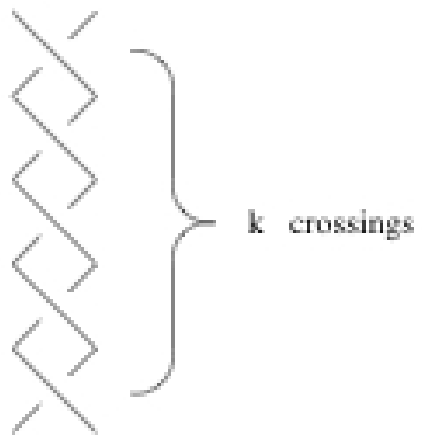}
\end{picture}
\end{center}

In this case $D_{-}$ is  reduced and alternating, so as in the first part $a_{k-3}(D_{-})=(i)^{j-1}$. As in the first case the signs work out to give
\begin{equation}
|a_{k-4}(D)-a_{k}(D)|= |a_{k-5}(D_{+})-a_{k-1}(D_{+})|+1
\end{equation}

Notice $D_{+}$ is  alternating and plus adequate, but it is no longer reduced. We get the following figure
\begin{center}
\begin{picture}(125,127)
\includegraphics{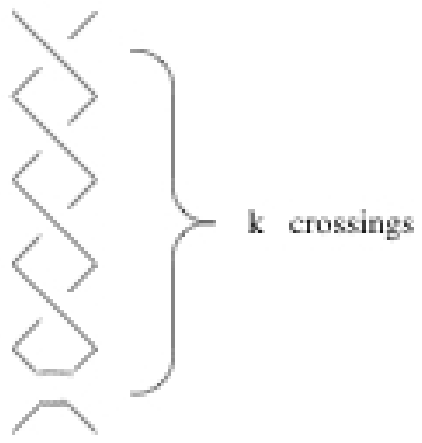}
\end{picture}
\end{center}

 Consider the coefficients of the lowest degree terms of the Kauffman bracket. The above lemma gives the following relation.
 
 \begin{equation}
|a_{l+4}(D)-a_{l}(D)|  =  |a_{l+3}(D_{+})-a_{l+5}(D_{-})+a_{l+1}(D_{-})| 
\end{equation}
It will be shown that $a_{l+3}(D_{+})=0$. First consider the state in $D_{+}$ where all crossing are smoothed negatively. The (k-1)-twists then give the following picture
\begin{center}
\begin{picture}(125,127)
\includegraphics{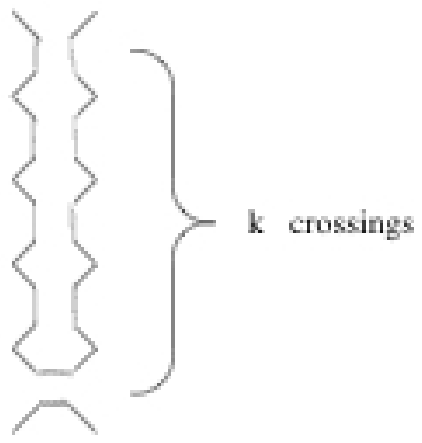}
\end{picture}
\end{center}
 
  The coefficient of the lowest degree term in the Kauffman bracket of $D_{+}$ is $a_{l+3}(D_{+})$. Consider the homology groups $H_{s,l+3}(D_{+})$ where s ranges between m+1 and j-1, by two. Recall, m is the lowest homological degree and j is the highest homological degree achieved in the homology of D. In fact m is the negative of the number of crossings and j is the number of crossings. The lowest and highest homological degrees achieved in the homology of $D_{-}$ are  m+1 and j-1 .
 
  The boundary operator preserves the polynomial degree of the chain complex. As such the overall complex $C(D_{+})$ decomposes into the direct sum of several complexes, one for each degree. According to a theorem of Hopf the Euler characteristic of the chain complex is the same as the Euler characteristic of it's homology groups. Thus it suffices to consider the complex $C_{l+3}(D_{+})$.
  
   Let S be an enhanced state with polynomial degree l+3 where all the circles are oriented counter clockwise. Let $S_{-}$ be the state of $ D_{+}$ where all crossing are smoothed negatively and all circles are oriented counterclockwise. Notice $S_{-}$ also has polynomial degree l+3. From the definition of polynomial degree the following holds for r where  $\frac{r}{2}$ is the number of negative smoothings changed to positive smoothings to arrive at the state S from $S_{-}$. 
 \begin{equation}
 l+3=\sigma(S_{-}) + 2 \tau(S_{-}) =(m+1)+2\tau(S_{-})=\sigma(S) + 2\tau(S)=((m+1)+r)+2\tau(S)
 \end{equation}

 Notice $\tau(S)$ and $\tau(S_{-})$ are actually the negative of the number of circles in each state respectively. From the above equation, $\tau(S_{-})-\tau(S)=\frac{r}{2}$.

 Suppose we did not take the state S with all circles oriented counterclockwise. Let $S'$ be the same smoothing as S, but with t circles oriented clockwise. Then $\tau(S')=\tau(S)+2t$. Therefore the polynomial degree of $S'$ is 

\begin{equation}
((m+1)+r) +2\tau(S')=((m+1)+r)+2\tau(S) +4t \neq l+3
\end{equation} 
 Thus the state $S'$ is not in the sub-complex  $C_{l+3}(D_{+})$.  Therefore all states belonging to $C_{l+3}(D_{+})$ have all their circles oriented counter clockwise. 
 
 Notice that $C_{m+1,l+3}(D_{+})$ is generated by the state where all crossings have been smoothed negatively in $D_{+}$ and all circles oriented counterclockwise. Suppose that at only the  crossings in $D_{+}$, which correspond to the equivalence class  containing v in $D$, can we increase the number of circles by changing a crossing. There are (k-1) ways to increase the number of circles by one, that is by changing the smoothing at each of the (k-1) crossings that belong to the equivalence class that contained v. Thus the rank of $C_{(m+1)+2,l+3}(D_{+})$ is (k-1). Then there are $(k-1)\choose 2$ ways of creating two circles, that is, by changing two of the smoothings corresponding to crossings in the twist. We can say in general that the rank of $C_{(m+1)+r,l+3}(D_{+})$ is$(k-1)\choose \frac{r}{2}$. The above is true so long as there is no crossing in another equivalence which acts like those that are equivalent to v in D.  

 Suppose that there was some other crossing in $D_{+}$ that is not in the equivalence class of v in D, say w, and that by smoothing the crossing at w positively in $D_{+}$ the number of circles go up when compared to the state of $D_{+}$ with all crossings smoothed negatively. Consider the state of $D_{+}$ that has all crossings smoothed negatively. Since the number of circles goes up it must be that w joins a circle in the smoothing of $D_{+}$ to itself. This can happen in two ways, inside the circle, or outside the circle. 
\begin{center}
\begin{picture}(172,107)
\includegraphics{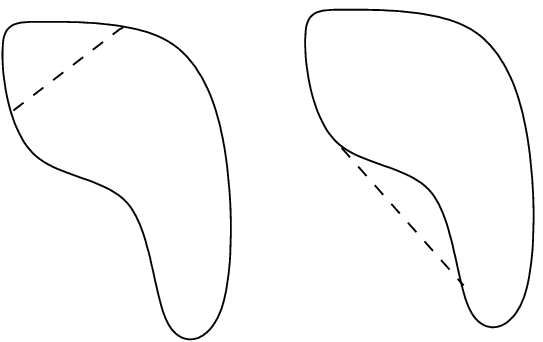}
\end{picture}
\end{center}
Suppose that the circle that w connects to itself is not a circle that involves v.
\begin{center}
\begin{picture}(280,188)
\includegraphics{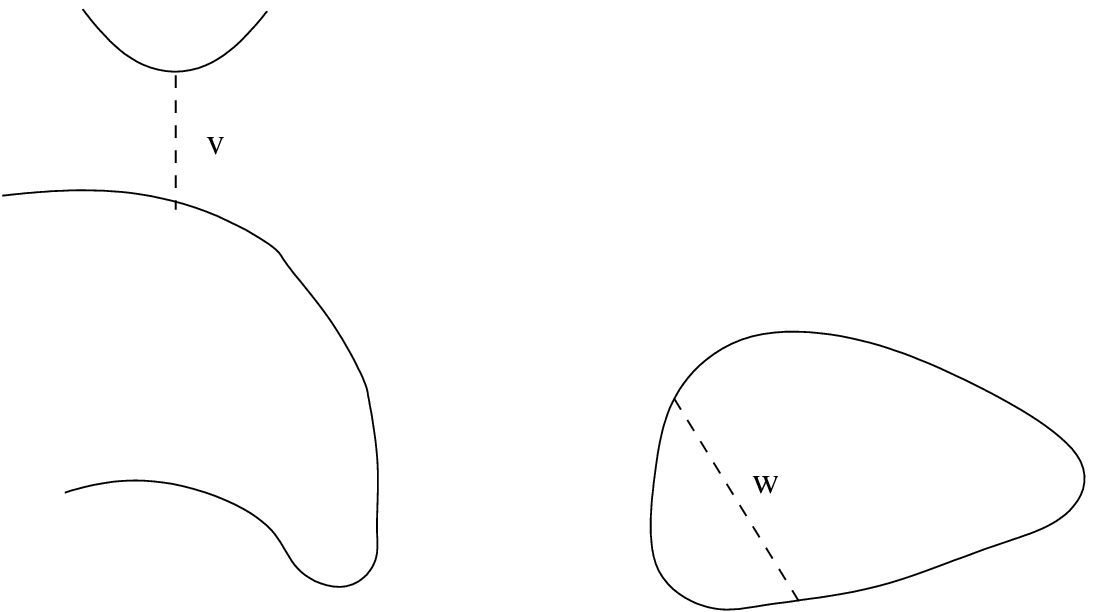}
\end{picture}
\end{center} 
 Then in the state of D with all crossings smoothed negatively, w would still connect a circle to itself. Thus if we smooth w positively the number of circles goes up. This is a contradiction to D being adequate.
Thus w connects a circle to itself which involves v . No crossing equivalent to w can connect the circle to itself  as in figure (a) below. Otherwise the diagram is not reduced. If there were a crossing that ensures the diagram is reduced, as in figure (c) the diagram is not alternating since v is the only crossing with a positive smoothing. 
\begin{center}
\begin{picture}(257,148)
\includegraphics{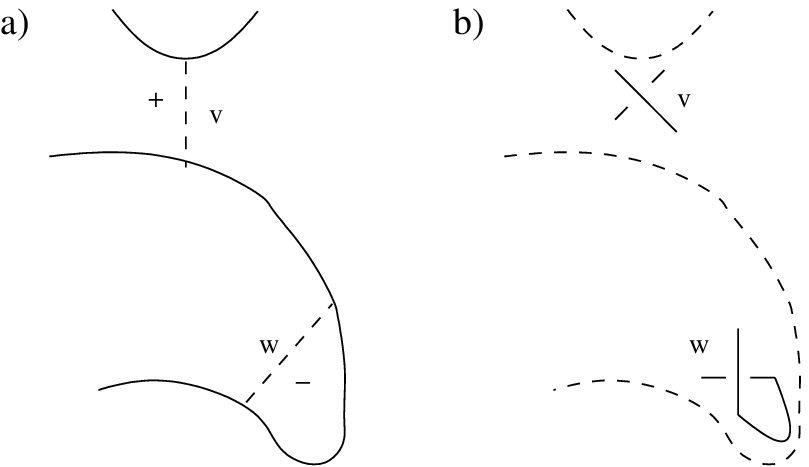}
\end{picture}
\end{center}

\begin{center}
\begin{picture}(273,186)
\includegraphics{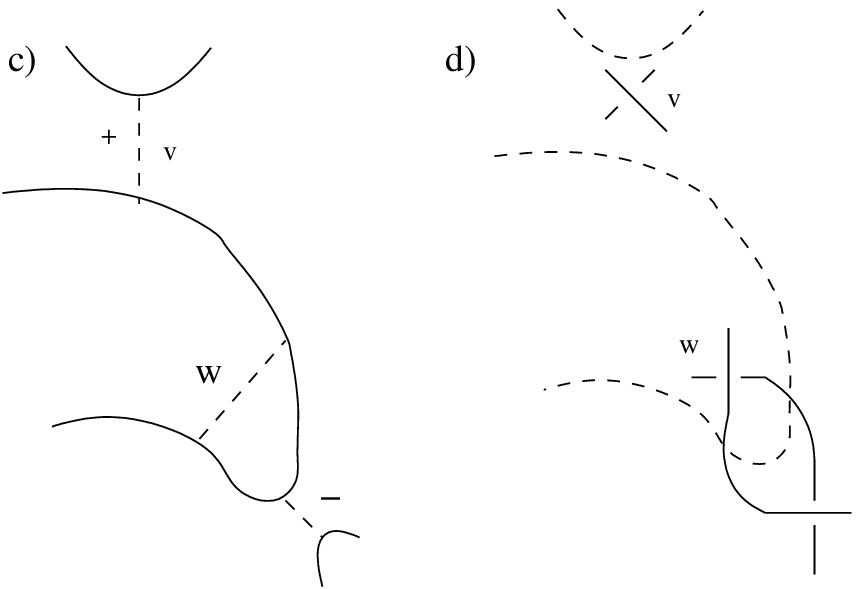}
\end{picture}
\end{center}
Thus the only way that we can see w in this circle while maintaining a reduced alternating diagram is as indicated in figure  below. Then v and w are equivalent contradicting our assumption.
 \begin{center}
\begin{picture}(268,157)
\includegraphics{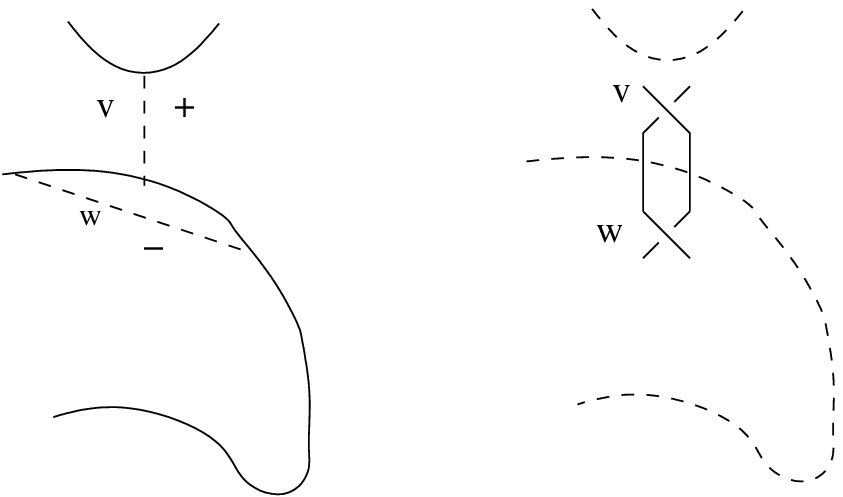}
\end{picture}
\end{center}
The same thing happens when the crossing is outside the circle. 
Thus only changing the smoothing at the k-1 crossings in the twist that contained v  can increase the number of circles.  Note that $\frac{r}{2}$ runs from  0 to (k-1). It follows that
\begin{equation}
a_{l+3}(D_{+})=\sum_{c=m+1...j-1 \ by 2 \ } (i)^{c} rank(H_{c,l+3}(D_{+})=(i)^{m}\sum_{c=0..k-1}(-1)^{c} {k-1\choose c} =0
\end{equation}
 
and

\begin{multline}
|a_{l+4}(D)-a_{l}(D)|=\\
|a_{l+3}(D_{+})-a_{l+5}(D_{-})+a_{l+1}(D_{-})| =|-a_{l+5}(D_{-})+a_{l+1}(D_{-})|\\
\end{multline}

A positive smoothing at v deletes an edge in $\tilde{G^{*}}$. A negative smoothing at v  deletes one edge in a set of parallel edges so that $\tilde{G}$ does not change.

By the inductive hypothesis 
\begin{eqnarray}
|a_{j-4}(D)-a_{j}(D)|&=\psi(\tilde{G^{*}_{+}})+1&=\psi(\tilde{G^{*}})\\
|a_{l+4}(D)-a_{l}(D)|&=\psi(\tilde{G_{-}})=\psi(\tilde{G})
\end{eqnarray}

The proof is similar when the crossing we start with is the mirror image of v. 
\end{proof}

\begin{cor}

In the case when L is an alternating knot then \cite{dasb} showed $(|\tilde{E}|+|V|+1)+ 
(|\tilde{E^{*}}| +|V^{*}| +1)$=twist number of L. Thus $|a_{k-4}-a_{k}|+|a_{l+4}-a_{l}|$=twist number of L.  

\end{cor}


\begin{thebibliography}{0000}
\bibitem{dasb}Dasbach, O. , Lin, X.S. A Volumish Theorem for the Jones Polynomial,arXiv:math.GT/04034488 v1 2004
\bibitem{lick} Lickorish, W.B. Raymond. An Introduction to knot theory, New York: Springer, c1997.
\bibitem{viro} Viro, Oleg. Remarks on the Definition of Khovanov Homology,arXiv:math.GT/0202199 v1 2002
\end{thebibliography}
\end{document}